\documentclass[reqno,12pt]{amsart}

\usepackage{amsmath}
\usepackage{amscd}
\usepackage{amssymb}
\usepackage{enumerate}
\usepackage{amsfonts}
\usepackage{graphicx}
\usepackage[all]{xy}
\usepackage{mathrsfs}

\DeclareFontEncoding{OT2}{}{}

\numberwithin{equation}{section}

\newcommand{\rar}[1]{\stackrel{#1}{\longrightarrow}}
\newcommand{\xrar}[1]{\xrightarrow{#1}}

\newcommand{\isom}{\rar{\simeq}}
\newcommand{\eqv}{\rar{\sim}}

\newcommand{\into}{\hookrightarrow}
\newcommand{\onto}{\twoheadrightarrow}

\newcommand{\al}{\alpha} \newcommand{\be}{\beta} 
\newcommand{\Ga}{\Gamma} \newcommand{\de}{\delta} 
\newcommand{\la}{\lambda} \newcommand{\La}{\Lambda}
\newcommand{\ze}{\zeta} 
\newcommand{\eps}{\epsilon} \newcommand{\sg}{\sigma}

\newcommand{\bC}{{\mathbb C}}

\newcommand{\bN}{{\mathbb N}}

\newcommand{\bZ}{{\mathbb Z}}

\newcommand{\cD}{{\mathcal D}}
\newcommand{\cE}{{\mathcal E}}

\newcommand{\cO}{{\mathcal O}}

\newcommand{\cR}{{\mathcal R}}

\newcommand{\cU}{{\mathcal U}}

\newcommand{\sE}{{\mathscr E}}

\newcommand{\sH}{{\mathscr H}}

\newcommand{\sL}{{\mathscr L}}
\newcommand{\sM}{{\mathscr M}}
\newcommand{\sN}{{\mathscr N}}
\newcommand{\sO}{{\mathscr O}}
\newcommand{\sP}{{\mathscr P}}

\newcommand{\fX}{{\mathfrak X}}

\newcommand{\fg}{{\mathfrak g}}

\newcommand{\fl}{{\mathfrak l}}

\newcommand{\fp}{{\mathfrak p}}

\newcommand{\Qb}{\overline{Q}}

\newcommand{\mt}{\widetilde{M}}
\newcommand{\xt}{{\widetilde{X}}}

\newcommand{\dd}{\partial}

\newcommand{\abs}[1]{\vert #1\vert}

\newcommand{\Ker}{\operatorname{Ker}}
\newcommand{\End}{\operatorname{End}}
\newcommand{\Hom}{\operatorname{Hom}}
\newcommand{\Ext}{\operatorname{Ext}}

\newcommand{\Spec}{\operatorname{Spec}}
\newcommand{\Proj}{\operatorname{Proj}}

\newcommand{\rk}{\operatorname{rk}}

\newcommand{\ev}{\operatorname{ev}}

\newcommand{\Rep}{\operatorname{Rep}}
\newcommand{\gr}{\operatorname{gr}^{\bullet}}

\newcommand{\tens}{\otimes}

\newcommand{\st}{\,\big\vert\,}

\newcommand{\sbr}{\smallbreak}
\newcommand{\mbr}{\medbreak}

\newcommand{\lMod}[1]{#1\!-\!\text{Mod}}
\newcommand{\lmod}[1]{#1\!-\!\text{mod}}
\newcommand{\rMod}[1]{\text{Mod}\!-\!#1}
\newcommand{\rmod}[1]{\text{mod}\!-\!#1}
\newcommand{\Grmod}[1]{#1\!-\!\text{Grmod}}
\newcommand{\grmod}[1]{#1\!-\!\text{grmod}}
\newcommand{\Tors}[1]{#1\!-\!\text{Tors}}
\newcommand{\tors}[1]{#1\!-\!\text{tors}}
\newcommand{\Qgr}[1]{#1\!-\!\text{Qgr}}
\newcommand{\qgr}[1]{#1\!-\!\text{qgr}}

\newtheorem{thm}{Theorem}[section]

\newtheorem{lem}[thm]{Lemma}
\newtheorem{sublem}[thm]{Sublemma}
\newtheorem{prop}[thm]{Proposition}

\theoremstyle{remark}
\newtheorem{rem}[thm]{Remark}

\newtheorem{defin}[thm]{Definition}

\newcommand{\cm}{Cohen-Macaulay\ }

\newcommand{\inv}[1]{^{G,\,#1}}

\newcommand{\lp}{\La_+}
\newcommand{\lpp}{\La_{++}}

\newcommand{\cst}{\bC^\times}

\newcommand{\muz}{\mu^{-1}(0)}
\newcommand{\muzs}{\mu^{-1}(0)^{ss}_\chi}

\newcommand{\ten}[1]{^{\tens #1}}

\title[Quantization of minimal resolutions]{Quantization of minimal resolutions
of \\ Kleinian singularities}

\author[M.~Boyarchenko]{Mitya Boyarchenko}

\thanks{The author's research was partially supported by NSF grant
DMS-0401164. \\ {\em Address:} Department of Mathematics,
University of Chicago, Chicago, IL 60637 \\ {\em E-mail:} {\tt
mitya@math.uchicago.edu}}

\begin{document}

\begin{abstract}
In this paper we prove an analogue of a recent result of Gordon
and Stafford that relates the representation theory of certain
noncommutative deformations of the coordinate ring of the $n$-th
symmetric power of $\bC^2$ with the geometry of the Hilbert scheme
of $n$ points in $\bC^2$ through the formalism of $\bZ$-algebras.
Our work produces, for every {\em regular} noncommutative
deformation $\cO^\la$ of a Kleinian singularity $X=\bC^2/\Ga$, as
defined by Crawley-Boevey and Holland, a filtered $\bZ$-algebra
which is Morita equivalent to $\cO^\la$, such that the associated
graded $\bZ$-algebra is Morita equivalent to the minimal
resolution of $X$. The construction uses the description of the
algebras $\cO^\la$ as quantum Hamiltonian reductions, due to
Holland, and a GIT construction of minimal resolutions of $X$, due
to Cassens and Slodowy.
\end{abstract}

\maketitle

\setcounter{tocdepth}{1}

\tableofcontents


\section{Introduction}\label{s:intro}

\subsection{}\label{ss:O-lambda}
Let $\Ga\subset SL_2(\bC)$ be a finite subgroup, and
$\bigl\{\cO^\la\bigr\}$ the family of noncommutative deformations
of the singularity $X=\bC^2/\Ga$ constructed by Crawley-Boevey and
Holland \cite{CBH}. If $(Q,I)$ denotes the McKay quiver associated
to $\Ga$, then the parameter space for these deformations is
naturally identified with $\bC^I$. Explicitly, given
$\la\in\bC^I$, recall (\cite{CBH}, p.~606) that the {\em deformed
preprojective algebra} of $Q$ with parameter $\la$ is defined by
\[
\Pi^\la = \Pi^\la(Q) = \bC\Qb \biggl/ \Bigl\langle\Bigl\langle
\sum_{a\in Q} [a,a^*]-\la \Bigr\rangle\Bigr\rangle,
\]
and $\cO^\la$ is the spherical subalgebra $\cO^\la = e_0 \Pi^\la
e_0$. Here $\Qb$ denotes the {\em double of $Q$}, i.e., the quiver
obtained from $Q$ by adding an arrow $a^*$ for each arrow $a\in Q$
such that the tail (resp., head) of $a^*$ is the head (resp.,
tail) of $a$. We write $\bC\Qb$ for the path algebra of $\Qb$, and
$e_i\in\bC\Qb$ for the idempotent corresponding to the vertex
$i\in I$; the extending vertex of $Q$ is denoted by $0\in I$.
Finally, $\la$ is identified with the element $\sum_{i\in I} \la_i
e_i\in\bC\Qb$.

\medbreak

We define a grading on the algebra $\bC\Qb$ by assigning degree
$0$ to each idempotent $e_i$, and degree $1$ to each arrow $a\in
Q$ and its opposite arrow $a^*$; in other words, the grading is by
the length of paths. The algebras $\Pi^\la$ and $\cO^\la$ inherit
natural filtrations from this grading. This filtration will be
used throughout the paper without explicit mention.

\subsection{}\label{ss:roots}
We let $\de\in\bN^I$ denote the minimal positive
imaginary root for $Q$, and we write
\[
\La=\bigl\{ \xi\in\bZ^I \st \xi\cdot\de=0 \bigr\},
\]
\[
\La_+ = \bigl\{ \xi\in\La \st \xi\cdot\al\geq 0 \text{ for every
positive Dynkin root } \al \bigr\},
\]
\[
\La_{++} = \bigl\{ \xi\in\La \st \xi\cdot\al > 0 \text{ for every
positive Dynkin root } \al \bigr\}.
\]
Hereafter, a {\em root} is an element of the root system
associated to the quiver $Q$, which in this case can be defined as
the set of all $\al\in\bZ^I\setminus\{0\}$ such that $q(\al)\leq
1$, where $q$ is the Tits form corresponding to $Q$. A root $\al$
is {\em Dynkin} if $\al\cdot\eps_0=0$, where $\eps_0\in\bZ^I$ is
the standard coordinate vector corresponding to the extending
vertex (in other words, the coordinate of $\al$ corresponding to
the extending vertex is zero). A root $\al$ is {\em real} (resp.,
{\em imaginary}) if $q(\al)=1$ (resp., $q(\al)=0$); note that a
Dynkin root is automatically real. There is a natural
identification of $\La$ with the weight lattice of the finite root
system associated to the Dynkin diagram obtained by deleting the
extending vertex, so that $\La_+$ (resp., $\La_{++}$) corresponds
to the set of dominant (resp., dominant regular) weights.

\subsection{} We denote by $\Rep(Q,\de)$ (resp., $\Rep(\Qb,\de)$)
the affine space of all representations of $Q$ (resp., $\Qb$) with
dimension vector $\de$. Using the trace pairing, $\Rep(\Qb,\de)$
is naturally identified with the cotangent bundle
$T^*\Rep(Q,\de)$. Let
\[
G=PGL(\de)=\left(\prod_{i\in I} GL(\de_i,\bC)\right)\biggl/\cst,
\]
where $\cst$ is embedded diagonally into the product. This is a
reductive algebraic group acting by conjugation on the varieties
$\Rep(Q,\de)$ and $\Rep(\Qb,\de)$, and we write
\[
\mu : \Rep(\Qb,\de) \rar{} \fg^*
\]
for the moment map for the action of $G$ on $\Rep(Q,\de)$ (see
\cite{CBH}, p.~606), where
\[
\fg=\operatorname{Lie}(G)=\fp\fg\fl(\de)=\left(\prod_{i\in I}
\fg\fl(\de_i,\bC)\right)\biggl/\bC
\]
is the Lie algebra of $G$.

\subsection{}\label{ss:Cass-Slod} Using the determinant maps
$\det:GL(\de_i,\bC)\to\cst$, we identify $\La$ with the group of
$1$-dimensional characters of $G$. Given $\chi\in\La_{++}$,
Cassens and Slodowy \cite{slod} construct a minimal resolution of
the Kleinian singularity $X$ as the projective morphism
\[
\xt:=\Proj S \rar{} \Spec S_0 \cong X,
\]
where $S$ is the graded algebra
\begin{equation}\label{e:cass-slod}
S=\bigoplus_{n\geq 0} S_n, \qquad
S_n=\bC[\mu^{-1}(0)]\inv{\chi^n}.
\end{equation}
Here $\bC[\mu^{-1}(0)]$ stands for the algebra of regular
functions on the {\em scheme-theoretic} fiber of $\mu$ over
$0\in\fg^*$ and $\bC[\muz]\inv{\chi^n}$ denotes the $G$-eigenspace
corresponding to the character $\chi^n$. Note that each component
$S_n$ of $S$ is itself graded, where the grading is induced by the
grading on $\bC[\mu^{-1}(0)]$, which in turn is induced by the
grading of $\bC[\Rep(\Qb,\de)]$ by the degree of polynomials (we
use the fact that $\Rep(\Qb,\de)$ is an affine space).

\medbreak

The minimal resolution $\xt\to X$ is studied in more detail in
Section \ref{s:minres}, where we also prove a result (Theorem
\ref{t:geom}) on the structure of the ring $S$ that, to the best
of our knowledge, does not appear in the existing literature.

\subsection{}\label{ss:b-la-chi} Let $\la\in\bC^I$ be such that $\la\cdot\de=1$ and
$\cO^\la$ is regular, i.e., has finite global dimension. By
Theorem 0.4 in \cite{CBH}, regularity in this case is equivalent
to either of the following two properties: $\cO^\la$ is Morita
equivalent to the deformed preprojective algebra $\Pi^\la$, or
$\la\cdot\al\neq 0$ for every Dynkin root $\al$. Using Holland's
description \cite{hol} of the algebras $\cO^\la$ as quantum
Hamiltonian reductions of the algebra of polynomial differential
operators on $\Rep(Q,\de)$, which is recalled in Section
\ref{s:holland}, we construct, for a given $\chi\in\La_{++}$, a
filtered $\bZ$-algebra $B^\la(\chi)$ such that
\begin{enumerate}[(1)]
\item $\lmod{\cO^\la}$ is naturally equivalent to
$\qgr{B^\la(\chi)}$, in a way compatible with filtrations, and
\item there is a natural isomorphism $\gr B^\la(\chi)\cong
\widehat{S}$, where $\widehat{S}$ is the $\bZ$-algebra associated
to the graded algebra $S$ defined by \eqref{e:cass-slod}.
\end{enumerate}
The construction is described in Section \ref{s:quantminres},
where we also state our main result, Theorem \ref{t:main}, which
is then proved in Section \ref{s:proof}.

\subsection{} The $\bZ$-algebra $B^\la(\chi)$ has the property that
$B^\la(\chi)_{0,0}=\cO^{\la+\xi}$ for a ``sufficiently large''
$\xi\in\lpp$, and then
\[
B^\la(\chi)_{n,n}=\cO^{\la+n\cdot\chi+\xi} \qquad \text{for all }
n\geq 0.
\]
Thus, the construction of $B^\la(\chi)$ depends on the choice of
$\xi$. However, this dependence is not very serious, since two
different choices of $\xi$ lead to naturally Morita equivalent
$\bZ$-algebras, which is why $\xi$ is omitted from the notation.
Furthermore, if $\cO^\la$ has no nonzero finite dimensional
modules, one can take $\xi=0$, and we conjecture that one can
always take $\xi=0$ as long as $\la$ is {\em dominant} (see Remark
\ref{r:xi}).

\medbreak

More importantly, $B^\la(\chi)$ also depends on the choice of
$\chi\in\lpp$. For this reason our quantization of the minimal
resolution $\xt$ may be called ``non-canonical.'' A more canonical
version of the quantization would consist of replacing
$B^\la(\chi)$ by a ``lower-triangular $\La$-algebra,'' bigraded by
$\lp$ instead of $\bZ_+$, and realizing $\xt$ as a suitable
``multi-$\Proj$'' of the $\lp$-graded ring $\oplus_{\chi\in\lp}
\bC[\muz]\inv{\chi}$. Note, however, that such a construction
cannot be obtained by a straightforward modification of the
results of the present paper. The most apparent reason for this is
that we repeatedly make crucial use of the following simple fact:
given a natural number $N$, every integer $n\geq 2N-1$ can be
written as a sum of integers that lie between $N$ and $2N-1$.
However, this fact has no suitable analogue for lattices other
than $\bZ$, in the sense that if $\rk\La\geq 2$, then $\lpp$ is
{\em not} finitely generated as a monoid. This issue will be
addressed in \cite{ii}.

\subsection{} A canonical (in the sense explained above)
quantization of minimal resolutions of Kleinian singularities of
type $A$ has been recently constructed by I.~Musson \cite{mus}.
His approach is very different from ours in that instead of using
Holland's results, he constructs a filtered $\bZ$-algebra that
deforms the minimal resolution by using the explicit description
of the latter as a toric variety (which replaces Cassens and
Slodowy's construction). In particular, this approach does not
generalize to other types of Kleinian singularities. Apart from
the basic theory of $\bZ$-algebras, our papers are completely
disjoint, and can be read independently.

\subsection{} We end the introduction by collecting some of the
notation and terminology that is used throughout the paper. A ring
is said to be {\em regular} if it has finite global homological
dimension. If $A$ is a ring, then $\lMod{A}$ and $\lmod{A}$ denote
the category of left (resp., left finitely generated) $A$-modules.
More generally, we follow the conventions of \cite{GS} and
\cite{mus}, in that any notation involving an uppercase letter
denotes a category of all modules with a certain property, whereas
the corresponding notation involving only lowercase letters
denotes the full subcategory consisting of finitely generated
modules.

\medbreak

In particular, if $A=\oplus_{n\geq 0} A_n$ is a graded algebra, we
denote by $\Grmod{A}$ the category of $\bZ$-graded (not
necessarily positively graded) $A$-modules. A graded $A$-module
$M$ is said to be {\em bounded} if $M_n=(0)$ for $n\gg 0$, and we
denote by $\Tors{A}\subseteq\Grmod{A}$ the full subcategory
consisting of modules that are unions of their bounded submodules.
The corresponding Serre quotient
\[
\Qgr{A} = \dfrac{\Grmod{A}}{\Tors{A}}
\]
plays the role of the ``category of quasi-coherent sheaves on
$\Proj A$.'' According to the previous paragraph, we have the
corresponding abelian categories $\grmod{A}$, $\tors{A}$ and
$\qgr{A}$. The notion of a (lower-triangular) {\em $\bZ$-algebra}
$B$, as well as the associated abelian categories $\Qgr{B}$,
$\qgr{B}$, etc., are introduced in Section \ref{s:morita}.

\mbr

An {\em algebra} will always mean for us an algebra over $\bC$,
and if $A,B$ are algebras, then an $(A,B)$-bimodule $M$ is
required to satisfy the condition that the two induced actions of
$\bC$ on $M$ coincide. All tensor products, unless specified
otherwise, will be taken over $\bC$. With the exception of
$\bZ$-algebras, all rings are assumed to have a multiplicative
identity, and all modules are assumed to be unital.

\subsection{Acknowledgements}
I am greatly indebted to Victor Ginzburg for introducing me to
this area of research, stating the problem, and constant
encouragement and attention to my work. During our numerous
conversations he suggested several ideas that were crucial for
some of the proofs appearing in this paper. I am also grateful to
Dennis Gaitsgory for motivating discussions during the early
stages of this work. Finally, I would like to thank Ian Musson for
explaining his work to me and sharing his preprint \cite{mus}
before it was made available to the general audience.

\sbr

I am thankful to Victor Ginzburg, Iain Gordon and Ian Musson for
making helpful comments, suggesting improvements, and pointing out
several misprints in the first version of the paper.


\section{A study of the minimal resolution}\label{s:minres}

\subsection{} Recall the notation $X=\bC^2/\Ga$
and $\mu:\Rep(\Qb,\de)\to\fg^*$ defined in the introduction. We
will denote by $\mu^{-1}(0)$ the scheme-theoretic fiber of $\mu$
over $0$. Thus
\begin{equation}\label{e:zerofiber}
\bC[\mu^{-1}(0)] =
\frac{\bC[\Rep(\Qb,\de)]}{\bC[\Rep(\Qb,\de)]\cdot\fg},
\end{equation}
where, by abuse of notation, $\fg$ denotes the linear subspace of
$\bC[\Rep(\Qb,\de)]$ obtained by pulling back via $\mu$ the
elements of $\fg$ viewed as linear functions on $\fg^*$. However,
by a result of Crawley-Boevey \cite{cb}, the scheme $\mu^{-1}(0)$
is in fact reduced and irreducible. (The reason for defining
$\mu^{-1}(0)$ as the scheme-theoretic fiber is that
\eqref{e:zerofiber} will be important for us later on.) We define
\[
R = \bC[\mu^{-1}(0)]^G.
\]
It is well known that
\begin{equation}\label{e:21}
\Spec R = \Spec \bC[\muz]^G = \mu^{-1}(0)//G \cong X.
\end{equation}
In particular, $R$ is a normal, $2$-dimensional, commutative
Gorenstein domain.

\subsection{}\label{ss:cass-slod} Cassens and Slodowy (\cite{slod}, \S7)
explain that a minimal resolution of $X$ can be constructed as a
GIT quotient
\begin{equation}\label{e:22}
\xt = \muzs/G
\end{equation}
for any $\chi\in\lpp$, where $\muzs$ denotes the open subset of
$\muz$ consisting of the points {\em semistable} with respect to
$\chi$. Moreover, they prove that for each such $\chi$,
\begin{enumerate}[(A)]
\item $\muzs=\muz^s_\chi$, the set of {\em stable} points with
respect to $\chi$, and
\item the action of $G$ on $\muz^s_\chi$ is {\em free} (recall
that, a priori, the action of $G$ on the set of stable points only
needs to have finite stabilizers; in our situation, however, all
stabilizers turn out to be trivial).
\end{enumerate}

Moreover, \eqref{e:21} and \eqref{e:22} lead to the description of
the resolution $\xt\to X$ as the natural map
\[
\Proj S \rar{} \Spec S_0,
\]
where $S$ is the graded ring
\[
S = \bigoplus_{n\geq 0} T\inv{\chi^n},
\]
where $T=\bC[\muz]$.

\begin{lem}\label{l:fingen}
The ring $S$ is finitely generated.
\end{lem}
\begin{proof}
This is a very general statement. Write $Y=\Spec T$, and consider
the induced action of $G$ on $Y\times\bC$, where the action on the
first factor is the given one, and the action on $\bC$ is via
$\chi$. The ring $\bC[Y\times\bC]=T\otimes\bC[z]$ has the obvious
grading by the degree of polynomials with respect to $z$, and we
clearly have an isomorphism of graded algebras
\[
\bC[Y\times\bC]^G \cong \bigoplus_{n\geq 0} T^{G,\chi^n}.
\]
In particular, the algebra on the RHS is finitely generated (here
we have used the fact that $G$ is reductive).
\end{proof}

\subsection{}\label{ss:geom-thm} The main goal of this section is to obtain some more
detailed information on the ring $S$, in the form of the following
result.
\begin{thm}\label{t:geom}
Let $p:\mu^{-1}(0)^{ss}_\chi \to\xt$ denote the quotient map.
\begin{enumerate}[(1)]
\item There exists a unique line bundle $\sL$ on $\xt$ such that
$p^*\sL$ is the trivial line bundle on $\mu^{-1}(0)^{ss}_\chi$
equipped with the $G$-linearization given by the character $\chi$.
Moreover, $\sL$ is ample.
\item The induced map
\[
S_n=\bC[\muz]\inv{\chi^n}\rar{} \Ga(\xt,\sL^{\tens n})
\]
is an isomorphism for sufficiently large $n$. In particular, $S_n$
is a torsion-free $S_0$-module of generic rank $1$ for
sufficiently large $n$.
\item The multiplication map
\[
S_m\tens S_n \rar{} S_{m+n}
\]
is surjective for sufficiently large $m$ and $n$.
\end{enumerate}
\end{thm}

\subsection{} It will be clear from the proof of the theorem that
essentially the only properties that we use are the fact that $G$
is reductive, statements (A) and (B) of \S\ref{ss:cass-slod}, and
the fact that $\muzs$ is dense in $\muz$. Thus the theorem could
be stated and proved in a much more general context, where $\muz$
is replaced by any affine variety $Y$ with an action of $G$
satisfying properties (A) and (B) of \S\ref{ss:cass-slod}, such
that the set of semistable points $Y^{ss}_\chi$ is dense in $Y$.

\subsection{}\label{ss:L} We begin the proof of Theorem \ref{t:geom} by
observing that (A) and (B) imply that the quotient map $p$ is a
principal $G$-bundle. Now it is easy to see that the notion of a
$G$-linearization for a coherent sheaf on $\muzs$ is equivalent to
the notion of a descent datum for the (flat) morphism $p$. Hence
the first statement of part (1) of the theorem follows immediately
from flat descent theory. Another consequence of descent theory is
that for any line bundle $\sM$ on $\xt$, we have
\begin{equation}\label{e:23}
\Ga(\xt,\sM) = \Ga(\muzs,p^*\sM)^G.
\end{equation}
Indeed, the LHS of \eqref{e:23} coincides with
$\Hom(\sO_{\xt},\sM)$. But $p^*\sO_{\xt}$ is the trivial line
bundle on $\muzs$ equipped with the trivial $G$-linearization, and
descent theory for morphisms implies that
\[
\Hom(\sO_\xt,\sM) = \Hom_{G-equiv}(p^*\sO_\xt,p^*\sM),
\]
which proves \eqref{e:23}.

\subsection{} The construction of the quotient \eqref{e:22} given
in \cite{git} shows that for some $N\in\bN$, there exists an ample
line bundle $\sL'$ on $\xt$ such that $p^*\sL'$ is the trivial
line bundle on $\muzs$ equipped with the $G$-linearization given
by $\chi^N$. Moreover, $\sL'=\sO_\xt(N)$ for the description of
$\xt$ as $\Proj S$. Now the discussion in \S\ref{ss:L} implies
that $\sL'=\sL^{\tens N}$; in particular, $\sL$ itself is ample,
completing the proof of part (1) of the theorem.

\medbreak

Replacing $N$ by one of its multiples if necessary, we may assume
that
\[
S_{jN} = (S_N)^j \qquad \text{for all } j\geq 1;
\]
this follows from the fact that $S$ is finitely generated (Lemma
\ref{l:fingen}) and Lemma 2.1.6(v) of \cite{ega2}. Similarly, we
may assume that $\sL^{\tens N}$ is very ample (using Proposition
4.5.10(ii) of loc. cit.). And, finally, we may assume that the
natural map
\[
S_{jN}\rar{} \Ga(\xt,\sL^{\tens jN})
\]
is an isomorphism for all $j\geq 1$. From now on we fix $N\in\bN$
satisfying all the properties listed above.

\subsection{} The line bundle $\sL^{\tens n}$ is very ample for
all sufficiently large $n$ (loc. cit., Proposition 4.5.10(ii)). In
particular, there exists $d_1\in\bN$ such that each of the bundles
\begin{equation}\label{e:24}
\sL^{\tens d_1 N}, \sL\ten{(d_1 N +1)}, \dotsc, \sL\ten{(d_1 N + N
-1)}
\end{equation}
is generated by global sections. But for any $n\in\bN$, we have,
from \eqref{e:23},
\[
\Ga(\xt,\sL\ten{n}) = \Ga(\muzs,\sO)\inv{\chi^n}.
\]
Now recall that $\muzs$ is precisely the set of points of $\muz$
where at least one element of $\bC[\muz]\inv{\chi^N}=S_N$ does not
vanish. In particular, if $\sg\in\Ga(\xt,\sL\ten{n})$, then there
exist $j\in\bN$ and finitely many elements $f_1,\dotsc,f_r\in
S_{jN}$ such that $f_i\sg\in S_{n+jN}$ for each $i$, and the open
sets $\{f_i\neq 0\}$ cover all of $\muzs$.

\medbreak

Since we are dealing with finitely many line bundles \eqref{e:24},
we deduce that there exists $d_2\in\bN$ such that for every $0\leq
k\leq N-1$, the line bundle $\sL\ten{(d_1 N+d_2 N+k)}$ is
generated by finitely many sections coming from the elements of
$S_{d_1N+d_2N+k}$.

\medbreak

For every $0\leq k\leq N-1$, let us now choose a finite
dimensional subspace
\[
V_k\subseteq S_{d_1 N+d_2
N+k}\subseteq\Ga(\xt,\sL\ten{(d_1N+d_2N+k)})
\]
of sections which generate the line bundle
$\sL\ten{(d_1N+d_2N+k)}$. These sections determine a surjection of
coherent sheaves
\begin{equation}\label{e:25}
\phi_k : \sO_\xt\tens_\bC V_k \rar{} \sL\ten{(d_1N+d_2N+k)}.
\end{equation}
Let $\sN_k$ denote the kernel of this surjection; it is a coherent
sheaf on $\xt$. Since $\sL\ten{N}$ is very ample, there exists
$d_3\in\bN$ such that
\begin{equation}\label{e:26}
H^1(\xt,\sL\ten{jN}\tens_{\sO_\xt}\sN_k) = 0
\end{equation}
for every $j\geq d_3$ and every $0\leq k\leq N-1$.

\subsection{} We can now prove part (2) of Theorem \ref{t:geom}.
Namely, every integer $n\geq (d_1+d_2+d_3)\cdot N$ can be written
as $n=jN+d_1N+d_2N+k$ for some (uniquely determined) $j\geq d_3$
and $0\leq k\leq N-1$. We have a short exact sequence, induced by
\eqref{e:25}:
\[
0 \rar{} \sL\ten{jN}\tens_{\sO_\xt}\sN_k \rar{}
\sL\ten{jN}\tens_\bC V_k \rar{} \sL\ten{n} \rar{} 0.
\]
Applying the long exact cohomology sequence and using
\eqref{e:26}, we see that the map
\[
\Ga(\xt,\sL\ten{jN}\tens_\bC V_k) \rar{} \Ga(\xt,\sL\ten{n})
\]
is surjective. But
\[
\Ga(\xt,\sL\ten{jN}\tens_\bC V_k) = \Ga(\xt,\sL\ten{jN})\tens_\bC
V_k = S_{jN}\tens_\bC V_k \subseteq S_{jN}\tens_\bC
S_{d_1N+d_2N+k},
\]
and so, a fortiori, the natural map
\[
S_n \rar{} \Ga(\xt,\sL\ten{n})
\]
is surjective for all $n\geq (d_1+d_2+d_3)\cdot N$. Furthermore,
this map is injective for all $n$ because $\muzs$ is dense in
$\muz$ (since $\muz$ is irreducible). Observe moreover that
$\Ga(\xt,\sL\ten{n})$ is a finitely generated $\bC[\muz]$-module
of generic rank $1$, since $\xt\to X$ is a projective birational
map which is an isomorphism away from the fiber over the singular
point $0\in X$. This proves part (2) of the theorem.

\subsection{} Finally, the proof of part (3) of Theorem
\ref{t:geom} is a variation of the idea used above, with an extra
trick. It is best to formulate it as a general statement. In view
of parts (1) and (2), it is clear that part (3) follows from the
following
\begin{prop}\label{p:ample}
Let $Y$ be a scheme, projective over a (commutative) Noetherian
ring $A$, and let $\sL$ be an ample invertible sheaf on $Y$. Then
the natural map
\begin{equation}\label{e:ga-tens}
\Ga(Y,\sL^{\tens m})\tens_A\Ga(Y,\sL^{\tens n}) \rar{}
\Ga(Y,\sL^{\tens(m+n)})
\end{equation}
is surjective for all sufficiently large $m$ and $n$.
\end{prop}
\begin{proof}
As before, the sheaves $\sL^{\tens n}$ are very ample for all
$n\gg 0$. In particular, there exists $N\in\bN$ such that for
every $N\leq k\leq 2N-1$, the sheaf $\sL\ten{k}$ is generated by a
finitely generated $A$-submodule of global sections
$V_k\subseteq\Ga(Y,\sL\ten{k})$. Let us write $\phi_k:\sO_Y\tens_A
V_k\to\sL\ten{k}$ for the corresponding surjection, and $\sN_k$
for its kernel.

\medbreak

Next, let $N'\in\bN$ be such that
\[
H^1(Y,\sL\ten{n}\tens_{\sO_Y}\sN_k) = 0 \quad \text{for all }
n\geq N' \text{ and } N\leq k\leq 2N-1.
\]
We claim that if $n\geq N'$ and $m\geq 2N-1$, then the map
\eqref{e:ga-tens} is surjective. Indeed, if $m\geq 2N-1$ is an
integer, it can clearly be written as a sum of the form
$m=k_1+\dotsb+k_r$, where $r\geq 1$ and $N\leq k_j\leq 2N-1$ for
every $j$. Now for $0\leq j\leq r-1$, we have an exact sequence
\[
0 \rar{} \sL\ten{(n+k_1+\dotsb+k_j)}\tens_{\sO_Y}\sN_k \rar{}
\sL\ten{(n+k_1+\dotsb+k_j)}\tens_A V_{k_{j+1}} \rar{}
\sL\ten{(n+k_1+\dotsb+k_j+k_{j+1})} \rar{} 0,
\]
and since $n+k_1+\dotsb+k_j\geq n\geq N'$, we see that the induced
map
\[
\Ga\Bigl(Y,\sL\ten{(n+k_1+\dotsb+k_j)}\Bigr)\tens_A V_{k_{j+1}}
\rar{} \Ga\Bigl(Y,\sL\ten{(n+k_1+\dotsb+k_j+k_{j+1})}\Bigr)
\]
is surjective. A fortiori, this implies that the map
\[
\Ga\Bigl(Y,\sL\ten{(n+k_1+\dotsb+k_j)}\Bigr)\tens_A \Ga\Bigl(
Y,\sL\ten{k_{j+1}}\Bigr) \rar{}
\Ga\Bigl(Y,\sL\ten{(n+k_1+\dotsb+k_j+k_{j+1})}\Bigr)
\]
is surjective.

\medbreak

By induction, we deduce that the map
\[
\Ga(Y,\sL\ten{n}) \tens_A \Ga(Y,\sL\ten{k_1}) \tens_A \dotsb
\tens_A \Ga(Y,\sL\ten{k_r}) \rar{} \Ga(Y,\sL\ten{(m+n)})
\]
is surjective. A fortiori, this implies that the map
\eqref{e:ga-tens} is surjective, completing the proof of the
proposition (and thereby of Theorem \ref{t:geom}).
\end{proof}


\section{Quantization of Kleinian singularities}\label{s:holland}

\subsection{} In this section we recall some results of M.P.~Holland
\cite{hol} that are crucial for our construction of the
quantization of the minimal resolution. We will use the notation
\[\Ga,\ I,\ Q,\ \de,\ G,\ \fg,\ \mu,\text{ etc.}\] defined in the
introduction. If $a\in Q$ is an arrow, we write $t(a),h(a)\in I$
for the tail and head of $a$, respectively. The {\em defect}
$\dd\in\bZ^I$ is defined by
\[
\dd_i=-\de_i+\sum_{t(a)=i} \de_{h(a)} \qquad \text{for all } i\in
I.
\]
We identify $\bC^I_0:=\{\chi\in\bC^I\big\vert\chi\cdot\de=0\}$
with the space of $1$-dimensional characters of $\fg$ via the
various trace maps $\fg\fl(\de_i,\bC)\to\bC$. If $\chi\in\bC^I_0$,
we define a filtered algebra
\[
\cU^\chi = \frac{\cD\bigl(\Rep(Q,\de)\bigr)^G}{\Bigl[
\cD\bigl(\Rep(Q,\de)\bigr)\cdot (\iota-\chi)(\fg) \Bigr]^G},
\]
where $\cD\bigl(\Rep(Q,\de)\bigr)$ is the algebra of polynomial
differential operators on the affine space $\Rep(Q,\de)$, and
\[
\iota : \fg\rar{} \operatorname{Vect}\bigl( \Rep(Q,\de) \bigr)
\subset \cD\bigl(\Rep(Q,\de)\bigr)
\]
is the Lie algebra map induced by the $G$-action.

\medbreak

\noindent
\textbf{CAUTION.} For consistency with the filtration on the
algebras $\cO^\la$ introduced in \S\ref{ss:O-lambda}, we need to
use the Bernstein filtration on the algebra $\cD(\Rep(Q,\de))$
(instead of the more standard order filtration). Fortunately, as
remarked in \cite{hol}, the results of sections 2--4 of that paper
remain valid if the order filtration is replaced by the Bernstein
filtration. From now on it will be implicitly assumed that the
results of all constructions involving differential operators will
be equipped with filtrations induced from the Bernstein
filtration.

\begin{thm}[Holland]\label{t:hol1}
If $\la\in\bC^I$ is such that $\la\cdot\de=1$, then there is a
natural isomorphism of filtered algebras
\[
\cO^\la\cong \cU^{\la-\dd-\eps_0}
\]
where $\eps_0\in\bZ^I$ is the standard basis vector corresponding
to the extending vertex.
\end{thm}
\begin{proof}
See \cite{hol}, Corollary 4.7.
\end{proof}

\subsection{} The following result will also be important to us.
\begin{thm}\label{t:hol2}
If $\chi\in\bC^I_0$, there is a natural isomorphism
\[
\gr \left( \frac{\cD(\Rep(Q,\de))}{\cD(\Rep(Q,\de))\cdot
(\iota-\chi)(\fg)} \right) \cong
\frac{\gr\cD(\Rep(Q,\de))}{\bigl[\gr
\cD(\Rep(Q,\de))\bigr]\cdot\fg}
\]
of graded $\bC[T^*\Rep(Q,\de)]$-modules.
\end{thm}
\begin{proof}
Combine \cite{hol}, Proposition 2.4 with the fact that the moment
map $\mu:\Rep(\Qb,\de)\to\fg^*$ is flat (\cite{CBH}, Lemma 8.3).
\end{proof}


\section{Morita $\bZ$-algebras}\label{s:morita}

\subsection{} In this section we review the basic theory of
$\bZ$-algebras following \cite{GS}. We also give a detailed proof
of a strengthening of Lemma 5.5 of loc. cit. that is used in our
paper.
\begin{defin}\label{d:z-alg}
A {\em lower-triangular $\bZ$-algebra} is given by the following
collection of data.
\begin{enumerate}[(a)]
\item An abelian group $B$, bigraded by $\bZ$ in the following
way:
\[
B = \bigoplus_{i\geq j\geq 0} B_{ij}.
\]
\item An algebra structure on each of the additive groups
$B_i:=B_{ii}$.
\item A $(B_i,B_j)$-bimodule structure on $B_{ij}$ for each
$i>j\geq 0$.
\item A map of $(B_i,B_k)$-bimodules
\[
B_{ij}\otimes_{B_j} B_{jk} \rar{} B_{ik}
\]
for each $i>j>k\geq 0$.
\end{enumerate}
The multiplication on the additive group $B$ is defined so as to
imitate multiplication of lower-triangular matrices in the obvious
sense, and the maps in (d) are required to satisfy the obvious
compatibility conditions needed to make the multiplication on $B$
associative. Recall that the algebras $B_i$ are assumed to have
multiplicative identities; however, $B$ will almost never have a
multiplicative identity, since it is defined as an infinite direct
sum.
\end{defin}
Next we consider modules over $\bZ$-algebras.
\begin{defin}\label{d:z-alg-mod}
Let $B$ be a lower-triangular $\bZ$-algebra as in the definition
above. A {\em graded $B$-module} is given by the following
collection of data.
\begin{enumerate}[(a)]
\item A positively graded abelian group
\[
M = \bigoplus_{i\geq 0} M_i.
\]
\item A $B_i$-module structure on $M_i$ for each $i\geq 0$.
\item A map of left $B_i$-modules
\[
B_{ij}\otimes_{B_j} M_j \rar{} M_i
\]
for each $i>j\geq 0$.
\end{enumerate}
The action of $B$ on $M$ is defined so as to imitate
multiplication of column vectors by lower-triangular matrices on
the left, and the maps in (c) are required to satisfy the obvious
compatibility conditions needed to make $M$ into a left
$B$-module. Note also that for each $i\geq 0$, the element $1\in
B_i$ is assumed to act as the identity on $M_i$.
\end{defin}

\subsection{} With these definitions in hand, we construct the
categories $\Grmod{B}$, $\Tors{B}$, $\Qgr{B}$, $\grmod{B}$,
$\tors{B}$ and $\qgr{B}$ in the same way as for graded rings. If
$A=\oplus_{n\geq 0} A_n$ is a graded ring, we can associate to it
a lower-triangular $\bZ$-algebra $B=\widehat{A}$ by defining
$B_{ij}=A_{i-j}$ for $i\geq j\geq 0$. As explained in \S5.3 of
\cite{GS}, we then have natural equivalences of categories
\[
\Qgr{A}\rar{\sim}\Qgr{\widehat{A}} \quad \text{and} \quad
\qgr{A}\rar{\sim}\qgr{\widehat{A}}.
\]

\subsection{} We are now ready for the key definition; note that
it is {\em weaker} than the corresponding notion introduced in
\S5.4 of \cite{GS}.
\begin{defin}\label{d:morita}
A {\em Morita $\bZ$-algebra} is a lower-triangular $\bZ$-algebra
\[
B = \bigoplus_{i\geq j\geq 0} B_{ij}
\]
such that there exists $N\in\bN$ for which
\begin{enumerate}[(i)]
\item the $(B_i,B_j)$-bimodule $B_{ij}$ yields an equivalence
\[
\lmod{B_j}\rar{\sim} \lmod{B_i}
\]
whenever $i-j\geq N$; and
\item the multiplication map
\[
B_{ij}\otimes_{B_j} B_{jk} \rar{} B_{ik}
\]
is an isomorphism whenever $i-j,j-k\geq N$.
\end{enumerate}
\end{defin}

\begin{thm}\label{t:morita}
Suppose that $B$ is a Morita $\bZ$-algebra such that each $B_i$ is
a left Noetherian ring, and each $B_{ij}$ is a finitely generated
left $B_i$-module. Then
\begin{enumerate}[\rm (1)]
\item Each finitely generated graded left $B$-module is
graded-Noetherian.
\item The association
\[
\phi:M\longmapsto \bigoplus_{n\geq 0} B_{n,0}\otimes_{B_0} M
\]
induces an equivalence of categories
\[
\Phi:\lmod{B_0}\rar{\sim}\qgr{B}.
\]
\end{enumerate}
\end{thm}
\begin{proof}
We follow the proof of Lemma 5.5 in \cite{GS} closely.

\medbreak

\noindent
(1) Let $M$ be a finitely generated graded $B$-module. We have to
show that every graded submodule of $M$ is also finitely
generated. It is clear that $M$ is generated by finitely many
homogeneous elements, so it is enough to consider the case where
$M$ is generated by one homogeneous element, say of degree $a$. In
this case $M$ is a graded homomorphic image of $\oplus_{j\geq a}
B_{ja}$, so we assume, without loss of generality, that
\[
M = \bigoplus_{j\geq a} B_{ja}.
\]
Now let
\[
L = \bigoplus_{j\geq a} L_j \subseteq M
\]
be a graded submodule. We use the notation
\[
B_{ij}^* = \Hom_{\lmod{B_i}}(B_{ij},B_i),
\]
which is a $(B_j,B_i)$-bimodule. Let $N\in\bN$ be as in the
definition of a Morita $\bZ$-algebra. Then for $j\geq a+N$, we
have a chain of maps of left $B_a$-modules
\[
B_{ja}^*\tens_{B_j} L_j \into B_{ja}^*\tens_{B_j} M_j =
B_{ja}^*\tens_{B_j} B_{ja} \rar{\simeq} B_a,
\]
where the first map is injective because $B_{ja}^*$ is a
projective right $B_j$-module, and the second map is an
isomorphism by definition. We let
\[
X(j)\subseteq B_a
\]
be the image of the composition above. Since $B_a$ is assumed to
be Noetherian, there exists an integer $b\geq a+N$ such that
\[
\sum_{j\geq a+N} X(j) = \sum_{i=a+N}^b X(i) \subseteq B_a.
\]
Now for $k\geq a+N$, we have
\[
L_k\cong B_{ka}\tens_{B_a} B_{ka}^* \tens_{B_k} L_k,
\]
which means that
\[
L_k = B_{ka} X(k) \qquad \text{for } k\geq a+N,
\]
as submodules of $B_{ka}=M_k$. Thus, for $k\geq b+N$, we have
\begin{eqnarray*}
L_k &=& B_{ka} X(k) \subseteq \sum_{i=a+N}^b B_{ka} X(i) \\
&=& \sum_{i=a+N}^b B_{ki} B_{ia} X(i) = \sum_{i=a+N}^b B_{ki} L_i,
\end{eqnarray*}
where we have used the assumption that $B_{ki}\tens_{B_i} B_{ia}
\rar{\simeq} B_{ka}$ for $k-i,i-a\geq N$. Thus we see that
\[
\sum_{k\geq b+N} L_k \quad \text{is generated by} \quad
\sum_{i=a+N}^b L_i
\]
as a $B$-module. Finally, for $a\leq j\leq b+N$, $L_j$ is a
$B_j$-submodule of the finitely generated $B_j$-module
$M_j=B_{ja}$, and is therefore finitely generated, completing the
proof of (1).

\medbreak

\noindent
(2) Note that, in view of (1), the category $\tors{B}$ consists
precisely of the bounded finitely generated graded left
$B$-modules. This remark will be used without further mention
throughout the proof.

\smallbreak

We construct a functor
\[
\Theta : \qgr{B} \rar{} \lmod{B_0}
\]
as follows. Let $\mt\in\qgr{B}$ be the image of some
$M\in\grmod{B}$. Then $M$ is generated by $\oplus_{i=0}^a M_i$ as
a left $B$-module, for some $a\in\bN$. Thus
\[
M_j = \sum_{i=0}^a B_{ji} M_i \qquad \text{for each } j\geq a,
\]
and so, by assumption,
\[
M_k = \sum_{i=0}^a B_{ki} M_i = \sum_{i=0}^a B_{kj} B_{ji} M_i =
B_{kj} M_j
\]
whenever $j\geq a$ and $k\geq j+N$. So if $k-j,j-a\geq N$, we have
a chain of natural homomorphisms
\begin{eqnarray*}
B_{ja}^* \tens_{B_j} M_j && \isom B_{ja}^* \tens_{B_j} B_{kj}^*
\tens_{B_k} B_{kj} \tens_{B_j} M_j \\
&& \isom B_{ka}^* \tens_{B_k} \bigl( B_{kj}\tens_{B_j} M_j \bigr)
\onto B_{ka}^*\tens_{B_k} M_k.
\end{eqnarray*}
Let $\theta_{kj}$ denote the composition, which is a surjective
homomorphism of left $B_a$-modules:
\[
\theta_{kj} : B_{ja}^* \tens_{B_j} M_j \onto B_{ka}^* \tens_{B_k}
M_k.
\]
Since the construction is natural, it is clear that the
homomorphisms $\theta_{kj}$ satisfy the compatibility conditions
\[
\theta_{lj} = \theta_{lk}\circ\theta_{kj} \quad \text{for all }
l-k,k-j,j-a\geq N.
\]

\begin{sublem}\label{l:1}
There exists an integer $b\geq a+N$ such that the map
$\theta_{kj}$ is an isomorphism whenever $j\geq b$ and $k-j\geq
N$.
\end{sublem}

Suppose otherwise. Then there exist $j_1\geq a+N$ and $k_1\geq
j_1+N$ such that $\theta_{k_1 j_1}$ has nonzero kernel. Further,
there exist $j_2\geq k_1+N$ and $k_2\geq j_2+N$ such that
$\theta_{k_2 j_2}$ has nonzero kernel. Continuing this process
inductively, we see that there exists an infinite sequence of
surjective homomorphisms
\[
B_{j_1a}^*\tens_{B_{j_1}} M_{j_1} \xrar{\theta_{k_1j_1}}
B_{k_1a}^*\tens_{B_{k_1}} M_{k_1} \xrar{\theta_{j_2k_1}}
B_{j_2a}^*\tens_{B_{j_2}} M_{j_2} \xrar{\theta_{k_2j_2}} \dotsb
\]
that never stabilizes to a sequence of isomorphisms. But since
$B_{j_1a}^*$ is a finitely generated left $B_a$-module and
$M_{j_1}$ is a finitely generated left $B_{j_1}$-module, we see
that $B_{j_1a}^*\tens_{B_{j_1}} M_{j_1}$ is a finitely generated,
hence Noetherian, left $B_a$-module, which is a contradiction.

\medbreak

We now fix $b\in\bN$ satisfying the conditions of the sublemma. We
define
\[
L = B_{b0}^* \tens_{B_b} M_b\in\lmod{B_0}.
\]

\begin{sublem}\label{l:2}
We have $\Phi(L)\cong\mt$ in $\qgr{B}$, where the functor $\Phi$
is defined in the statement of the theorem.
\end{sublem}

Recall that $\Phi(L)$ is the image under the quotient map of
\[
\phi(L) = \bigoplus_{n\geq 0} B_{n0}\tens_{B_0} B_{b0}^*
\tens_{B_b} M_b.
\]
Note that for $n\geq b+N$, we have the natural isomorphism of
$(B_n,B_0)$-bimodules
\[
B_{nb}\tens_{B_b} B_{b0} \isom B_{n0}.
\]
Tensoring with $B_{b0}^*$ on the right, we obtain a natural
isomorphism
\[
B_{nb} \isom B_{n0}\tens_{B_0} B_{b0}^*.
\]
Similarly, we have
\[
B_{kj} \isom B_{ka}\tens_{B_a} B_{ja}^*
\]
for $j\geq b$, $k\geq j+N$. Thus we obtain natural isomorphisms
\[
M_k \isom B_{ka}\tens_{B_a} B_{ka}^*\tens_{B_k} M_k
\xrar{(1\tens\theta_{kj})^{-1}} B_{ka}\tens_{B_a}
B_{ja}^*\tens_{B_j} M_j \isom B_{kj}\tens_{B_j} M_j.
\]
In particular, taking $j=b$ and combining with the previous
remarks, we obtain natural isomorphisms
\[
B_{n0}\tens_{B_0} B_{b0}^* \tens_{B_b} M_b \isom B_{nb}\tens_{B_b}
M_b \isom M_n
\]
for all $n\geq b+N$. By the naturality of all the constructions,
we see that we have defined an isomorphism of graded left
$B$-modules
\[
\phi(L)_{\geq b+N} \isom M_{\geq b+N}=\bigoplus_{n\geq b+N} M_n,
\]
completing the proof of the sublemma.

\medbreak

We note that, by similar arguments, replacing $b$ with any $b'\geq
b$ yields a module $L'=B_{b'0}^*\tens_{B_{b'}}
M_{b'}\in\lmod{B_0}$ that is canonically isomorphic to $L$. In
particular, $L$ only depends on the image $\mt$ of $M$ in
$\qgr{B}$, up to isomorphism.

\medbreak

Going back to the construction of the functor
\[
\Theta : \qgr{B} \rar{} \lmod{B_0},
\]
we set $\Theta(\mt)=L$ with the notation above. We need to explain
how $\Theta$ acts on morphisms. Suppose $M,M'\in\grmod{B}$ and
$L=\Theta(\mt)$, $L'=\Theta(\mt')$. In the construction of $L$,
$L'$ given above, we pick $b\in\bN$ that works for both; thus
$L=B_{b0}^*\tens_{B_b} M_b$ and $L'=B_{b0}^*\tens_{B_b} M_b'$. Now
a morphism of graded $B$-modules
\[
f : M \rar{} M'
\]
gives in particular a $B_b$-module homomorphism $f_b:M_b\to M'_b$;
tensoring with $B_{b0}^*$ yields a $B_0$-module morphism
\[
\Theta(f):=1\tens f_b : L\rar{} L'.
\]
Finally, we need to check that the functors $\Phi$ and $\Theta$
are quasi-inverse to each other. It is essentially obvious that
$\Theta\circ\Phi$ is isomorphic to the identity functor on
$\lmod{B_0}$, so we only need to consider the composition
$\Phi\circ\Theta$.

\medbreak

In the proof of Sublemma \ref{l:2} we have constructed a natural
isomorphism
\[
\Phi(\Theta(\mt))\isom\mt
\]
for every $M\in\grmod{B}$. We only need to see that this
isomorphism is compatible with morphisms; the essential
observation is that if $f:M\to M'$ is a morphism in $\grmod{B}$
and $b$ is as in the previous paragraph, then $f_n$ is uniquely
determined by $f_b$ for all $n\geq b+N$, due to the commutativity
of the diagram
\[
\xymatrix{
B_{nb}\tens_{B_b} M_b \ar[d] \ar[rr]^{1\tens f_b} & \ \ & B_{nb}\tens_{B_b} M'_b\ar[d] \\
 M_n \ar[rr]^{f_n} & \ \  & M'_n }
\]
where the vertical arrows are isomorphisms. This completes the
proof of the theorem.
\end{proof}


\section{Quantization of the minimal
resolution}\label{s:quantminres}

\subsection{}\label{ss:bimods} In this section we use Holland's results described
in Section \ref{s:holland} to define the algebras $B^\la(\chi)$
mentioned in the introduction and state our main result. We use
the same notation as in the introduction and in Section
\ref{s:holland}. We let $\fX(G)$ denote the group of the algebraic
group homomorphisms $\zeta:G\to\cst$. Note that if
$\zeta\in\fX(G)$, its differential $d\zeta:\fg\to\bC$ can be
thought of as an element of $\bC^I_0$ (see Section
\ref{s:holland}).

\medbreak

Given $\zeta\in\fX(G)$ and $\chi\in\bC^I_0$, we define
\[
P_{\chi,\zeta}=\frac{\cD\bigl(\Rep(Q,\de)\bigr)\inv{\zeta}}{\Bigl[
\cD\bigl(\Rep(Q,\de)\bigr) \cdot (\iota-\chi)(\fg)
\Bigr]\inv{\ze}}.
\]
One immediately checks that the actions of $\cD(\Rep(Q,\de))^G$ on
$\cD(\Rep(Q,\de))\inv{\ze}$ by left and right multiplication
descend to a filtered $(\cU^{\chi+d\ze},\cU^\chi)$-bimodule
structure on $P_{\chi,\ze}$.

\subsection{}\label{ss:P-la-ze} If $\la\in\bC^I$ is such that $\la\cdot\de=1$, we
will write
\[
P^\la_\ze=P_{\la-\dd-\eps_0,\ze}.
\]
According to Theorem \ref{t:hol1}, we can think of $P^\la_\ze$ as
a filtered $(\cO^{\la+d\ze},\cO^\la)$-bimodule. Moreover, by
Theorem \ref{t:hol2}, we have an isomorphism of graded bimodules
\begin{equation}\label{e:assoc-gr}
\gr P^\la_\ze \cong \bC[\mu^{-1}(0)]\inv\ze.
\end{equation}
The following fact will be used implicitly in Section
\ref{s:proof}; it is needed in order to justify the use of
Proposition \ref{p:cohen-mac}(1).
\begin{lem}\label{l:good-filt}
The filtration on $P^\la_\ze$ induced by the Bernstein filtration
on differential operators is {\em good} in the sense of
\cite{bjork}, Definition 2.19.
\end{lem}
\begin{proof} In view of \eqref{e:assoc-gr} and the Remark after
the proof of Proposition 2.22 in \cite{bjork}, it suffices to show
that $\bC[\mu^{-1}(0)]\inv\ze$ is a finitely generated
$\bC[\mu^{-1}(0)]^G$-module for every $\ze$. Now the algebra
$\oplus_{n\geq 0} \bC[\muz]\inv{\ze^n}$ is finitely generated by
the argument given in the proof of Lemma \ref{l:fingen}, and
therefore our claim follows from Lemma 2.1.6(i) in \cite{ega2}.
\end{proof}

\subsection{}\label{ss:P-la-xi}
Observe now that differentiation of characters induces an
isomorphism of abelian groups
$d:\fX(G)\rar{\simeq}\La\subseteq\bC^I_0$; by abuse of notation,
if $\xi\in\La$, we will write
\[
P^\la_\xi=P^\la_\ze \quad\text{and}\quad \bC[\mu^{-1}(0)]\inv\xi =
\bC[\mu^{-1}(0)]\inv\ze
\]
where $\ze\in\fX(G)$ is such that $d\ze=\xi$. Now, given
$\chi\in\lpp$, we define a lower-triangular $\bZ$-algebra
$B(\la,\chi)$ by
\[
B(\la,\chi)_{ij} = \begin{cases} \cO^{\la+j\cdot\chi} & \text{ if
} i=j\geq 0, \\
P^{\la+j\cdot\chi}_{(i-j)\cdot\chi} & \text{ if } i>j\geq 0.
\end{cases}
\]
All the structure maps of this $\bZ$-algebra are induced by the
multiplication of elements of $\cD(\Rep(Q,\de))^G$, and all
compatibility conditions follow immediately from the associativity
of this multiplication. Observe that $B(\la,\chi)$ is naturally
filtered by the Bernstein filtration of differential operators.
(We leave the formulation of the general notion of a filtered
$\bZ$-algebra $B$ to the reader: each component $B_{ij}$ should be
positively filtered, and all the structure maps should be
compatible with the filtrations. See also \cite{gs2}.)

\subsection{}\label{ss:state-main-thm}
We are now ready to state our main result, whose proof occupies
Section \ref{s:proof}.
\begin{thm}\label{t:main} Let $\la\in\bC^I$ be such that
$\la\cdot\de=1$ and $\la\cdot\al\neq 0$ for every Dynkin root
$\al$ (i.e., the algebra $\cO^\la$ is regular). Given
$\chi\in\lpp$, there exists $\xi\in\lpp$ such that the
lower-triangular filtered $\bZ$-algebra
\[
B^\la(\chi):=B(\la+\xi,\chi),
\]
where $B(\la+\xi,\chi)$ is constructed above, satisfies the
properties (1) and (2) of \S\ref{ss:b-la-chi} of the introduction.
\end{thm}
This theorem answers a question posed by M.P.~Holland on page 831
of \cite{hol}. He constructs a quantization of a certain {\em
partial} resolution of the Kleinian singularity $X=\bC^2/\Ga$, and
asks for a quantization of the whole minimal resolution. It
appears that the formalism of $\bZ$-algebras provides an
appropriate framework for ``noncommutative projective geometry'',
and Theorem \ref{t:main} shows that the filtered $\bZ$-algebra
$B^\la(\chi)$ deserves to be called a noncommutative deformation
of the minimal resolution $\widetilde{X}\to X$.

\subsection{}\label{ss:gord-staf} On the other hand, as explained in \S1.3 of
\cite{GS}, Theorem \ref{t:main} naturally completes the following
diagram (which does not commute!):
\[
\xymatrix{
? \ar[d]_{\text{gr}} & \ \ & \lmod{\cO^\la} \ar[d]^{\text{gr}}\ar[ll]_{\sim} \\
 Coh(\widetilde{X}) & \ \  & Coh(X) \ar[ll]_{\text{pullback}} }
\]
Namely, the question mark can be replaced by the category
$\qgr{B^\la(\chi)}$. The corresponding fact for rational Cherednik
algebras of type $A$ has been proved in \cite{GS}, and then
applications of this result to the representation theory of these
algebras have been explored in \cite{gs2}. A similar study for the
algebras $\cO^\la$ will appear in \cite{ii}.


\section{Proof of the main theorem}\label{s:proof}

\subsection{} In this section we prove Theorem \ref{t:main}. First
we collect several unrelated facts used in the proof, some of
which may be well known to the experts in the corresponding areas
of mathematics, but which we explain here for the reader's
convenience (and for the lack of suitable references). The proof
of the theorem itself begins in \S\ref{ss:proof}.

\subsection{}\label{ss:geometric} We start with a general geometric statement.

\begin{prop}\label{p:torsfree}
Let $X$ be a normal affine irreducible algebraic surface over
$\bC$, and $\cE$ a torsion-free coherent sheaf on $X$. Then
$\Ext^n(\cE,\cO_X)$ is finite dimensional for all $n\geq 1$. In
addition, if $\cE$ has generic rank $1$, then
\begin{equation}\label{e:endos}
\End_{\cO_X}(\cE) = \Ga(X,\cO_X).
\end{equation}
\end{prop}
\begin{proof}
Note that $X$ has at worst finitely many singular points because
it is a normal surface.
As usual, we write $\cE^\vee=\sH om(\cE,\cO_X)$ for the dual
sheaf. The canonical map $\cE\to\cE^{\vee\vee}$ is known to be an
isomorphism away from the singular points of $X$ and possibly
another finite set of points; moreover, $\cE^{\vee\vee}$ is
locally free away from the singular points of $X$. This implies
that $\cE$ itself is also locally free away from a finite set of
points. Thus for $n\geq 1$, the sheaf $\sE xt^n(\cE,\cO_X)$ has
finite support. On the other hand, since $X$ is affine, the
category of coherent $\cO_X$-modules is equivalent to the category
of finitely generated $\bC[X]$-modules, and in particular
$\Ext^n(\cE,\cO_X)\cong\Ga(X,\sE xt^n(\cE,\cO_X))$. This implies
the first statement.


\medbreak

For the second statement, let $S\subset X$ be the finite set of
points where $\cE$ is not locally free. Then
\[
\Ga(X,\cO_X)\subseteq\End_{\cO_X}(\cE) \hookrightarrow
\End_{\cO_{X\setminus S}}\Bigl(\cE\big\vert_{X\setminus S}\Bigr)
=\Ga(X\setminus S,\cO_X) = \Ga(X,\cO_X),
\]
where the second inclusion follows from the assumption that $\cE$
is torsion-free, and the last equality follows from the assumption
that $X$ is normal. This completes the proof.
\end{proof}

\subsection{}\label{ss:o-lambda} We now prove two results concerning the algebras
$\cO^\la$ that can be deduced from \cite{CBH}.
\begin{prop}\label{p:findim}
Let $\la\in\bC^I$ be such that $\la\cdot\de=1$ and $\cO^\la$ is
regular. Given $d\in\bN$, the algebra $\cO^{\la+\xi}$ has no
nonzero modules of dimension $\leq d$ for all sufficiently large
$\xi\in\lp$.
\end{prop}
\begin{proof}
Obviously, it suffices to prove the statement above with the word
``nonzero'' replaced by the word ``simple''. We now recall from
\cite{CBH} that the deformed preprojective algebra $\Pi^\la$
defined in \S\ref{ss:O-lambda} is Morita equivalent to $\cO^\la$,
the equivalence being given explicitly by
\[
M\mapsto e_0\Pi^\la\otimes_{\Pi^\la} M = e_0 M,
\]
where $e_0\in\Pi^\la$ denotes the idempotent corresponding to the
extending vertex. Note that this functor acts as
$\al\mapsto\eps_0\cdot\al$ on the dimension vectors of the
modules.

\mbr

Now we know from \cite{CBH} that the dimension vectors of simple
finite dimensional $\Pi^\la$-modules are among the positive roots
$\al$ satisfying $\la\cdot\al=0$ (observe that all such roots are
necessarily real and non-Dynkin). Also, given $\xi\in\La$, the map
\[
\al\mapsto\al-(\xi\cdot\al)\de
\]
establishes a bijection of {\em finite} sets
\[
\Bigl\{ \text{roots } \al \text{ such that } \la\cdot\al=0 \Bigr\}
\rar{\sim} \Bigl\{ \text{roots } \be \text{ such that }
(\la+\xi)\cdot\be=0 \Bigr\}.
\]

\medbreak

If $\al$ is a real root, let us write
$\al'=\al-(\eps_0\cdot\al)\de$, which is a Dynkin root. Observe
that if $\xi\in\La$, then we have $\xi\cdot\al=\xi\cdot\al'$.
Thus, if $\be=\al-(\xi\cdot\al)\de$, then
\begin{equation}\label{e:albe}
\eps_0\cdot\be = \eps_0\cdot\al -
\xi\cdot\al'.
\end{equation}
Now, using finiteness, choose $N\in\bN$ such that
$\abs{\eps_0\cdot\al}\leq N$ for every root $\al$ with
$\la\cdot\al=0$, and choose $\xi_0\in\lp$ such that
\[
\xi_0\cdot\psi>N+d \quad \text{for every positive Dynkin root }
\psi.
\]
It follows that if $\xi\geq\xi_0$ and $\al,\be$ are related as
above, then \eqref{e:albe} implies
\[
\abs{\eps_0\cdot\be}\geq\abs{\xi\cdot\al'}-\abs{\eps\cdot\al}>N+d-N=d.
\]
Finally, note that $\cO^{\la+\xi}$ is regular for sufficiently
large $\xi$, since, according to \cite{CBH}, we only need to
choose $\xi$ large enough so that $(\la+\xi)\cdot\al\neq 0$ for
every Dynkin root $\al$. It then follows from the discussion above
that the possible dimensions of the simple finite dimensional
$\cO^{\la+\xi}$ modules are among the integers $\eps_0\cdot\be$,
where $\be$ is a positive root with $(\la+\xi)\cdot\be=0$, which
completes the proof of the proposition.
\end{proof}
\begin{rem}\label{r:findimright}
We have stated and proved the proposition above for left
$\cO^\la$-modules. However, the same result also holds for right
modules. Indeed, if $(Q,I)$ is any quiver, it is easy to see that
for any $\la\in\bC^I$ there is a natural isomorphism between
$\Pi^{-\la}(Q)$ and the opposite algebra of $\Pi^\la(Q)$ that
preserves the idempotents corresponding to the vertices of $Q$. In
particular, it follows from \cite{CBH} that the dimension vectors
of simple left $\Pi^\la(Q)$-modules and simple right
$\Pi^\la(Q)$-modules are the same.
\end{rem}

\begin{prop}\label{p:shifts} Let $\la\in\bC^I$ and $\xi\in\La$ be
such that $\la\cdot\de=1$ and the algebras $\cO^\la$ and
$\cO^{\la+\xi}$ are both regular. Then there exists a Morita
equivalence
\begin{equation}\label{e:shifts}
\lmod{\cO^\la}\rar{\sim}\lmod{\cO^{\la+\xi}},
\end{equation}
compatible with filtrations in the obvious sense.
\end{prop}
\begin{proof}
Under the assumptions of the proposition, the results of
\cite{CBH} imply that the functors $\Pi^\la e_0\tens_{\cO^\la}-\ $
and $e_0\Pi^{\la+\xi}\tens_{\Pi^{\la+\xi}}-\ $ provide equivalence
of categories
\[
\lmod{\cO^\la} \eqv \lmod{\Pi^\la} \quad \text{and} \quad
\lmod{\Pi^{\la+\xi}} \eqv \lmod{\cO^{\la+\xi}}
\]
that are compatible with filtrations. Hence, if we can show that
there is an equivalence $\lmod{\Pi^\la}\eqv\lmod{\Pi^{\la+\xi}}$
that is also compatible with filtrations, we can define
\eqref{e:shifts} as the composition of these three equivalences.

\mbr

Now let $E$ denote the affine space of all $\la\in\bC^I$ such that
$\la\cdot\de=1$. Each simple reflection $s_i$ corresponding to a
vertex $i\in I$ of $Q$ defines an automorphism $r_i:E\to E$, and
the reflection functors of \cite{CBH} provide an equivalence
$\lmod{\Pi^\la}\eqv\lmod{\Pi^{r_i\la}}$ for every $\la\in E$. It
it clear from the construction given in loc.~cit. that this
equivalence is compatible with filtrations. On the other hand, if
$\phi:Q\to Q$ is an automorphism of the underlying graph of $Q$,
it also induces an automorphism $\phi^*:E\to E$, and it is easy to
see that there is a natural isomorphism
$\Pi^\la\to\Pi^{\phi^*\la}$ of filtered algebras, for any $\la\in
E$. Thus we have reduced the proof of the proposition to Lemma
\ref{l:shifts-preproj} below.
\end{proof}

\begin{lem}\label{l:shifts-preproj}
Given $\xi\in\La$, the map $\la\mapsto\la+\xi$ can be written as a
composition of simple reflections and automorphisms of the graph
$Q$.
\end{lem}
\begin{proof}
We use some standard facts about root systems and affine Weyl
groups that can be found in Chapter VI of \cite{bour}. Consider
the vector space $V=(\bZ^I/\bZ\de)\tens_\bZ\bC$; it is well known
that the image of the set of real roots for $Q$ under the
projection map $\bZ^I\to\bZ^I/\bZ\de$ is a reduced root system
$\cR$ in the space $V$, in the sense of loc.~cit., Chapter VI,
\S1.4. Now $V^*$ is naturally identified with $\La\tens_\bZ\bC$,
and $E$ can be viewed as an affine space for the vector space
$V^*$. Let $W_{ext}$ denote the group of automorphisms of $E$
generated by the translations by the elements of $\La$ and by the
Weyl group $W_{fin}$ of the root system $\cR$; sometimes $W_{ext}$
is called the {\em extended Weyl group} of the root system $\cR$.
It follows from the results of loc.~cit., Chapter VI, \S2.1 and
\S2.3, that $W_{ext}$ has the alternate description as the group
of automorphisms of $E$ generated by the {\em affine} Weyl group
$W_{aff}$ (which by definition is generated by the simple
reflections corresponding to all vertices of $Q$; it is called
simply the Weyl group of $Q$ in \cite{CBH}) and the group of
automorphisms of the graph $Q$. This proves the lemma.
\end{proof}

\begin{rem}\label{r:shifts} As we will see below, our proof of
Theorem \ref{t:main} relies heavily on the fact that the ``shift
functors''
\begin{equation}\label{e:shift-functs}
P^\la_\xi \tens_{\cO^\la} - \ : \ \lmod{\cO^\la} \rar{}
\lmod{\cO^{\la+\xi}},
\end{equation}
where the bimodules $P^\la_\xi$ have been defined in
\S\ref{ss:P-la-xi}, are equivalences of categories for
``sufficiently large'' $\la$ and $\xi$. Even though the proof of
Proposition \ref{p:shifts} provides a definition of shift
functors, it seems impractical to use this result for quantization
of minimal resolutions Kleinian singularities, since it is hard to
compute explicitly the associated graded spaces of the bimodules
defining the equivalences \eqref{e:shifts}, and to control the
compositions of these equivalences. It is not known to us if the
shift functors of Proposition \ref{p:shifts} are isomorphic to the
shift functors \eqref{e:shift-functs}.
\end{rem}

\subsection{}\label{ss:morita} Next we prove a characterization of Morita
equivalence which is somewhat different from the ones found in
standard textbooks, and for which we don't know a reference.
\begin{prop}\label{p:morita}
Let $A$, $B$ be rings, and let $P$ be an $(A,B)$-bimodule such
that the natural ring homomorphisms
\[
A\rar{} \End_{\rMod{B}}(P) \quad \text{and} \quad B\rar{}
\End_{\lMod{A}}(P)^{op}
\]
are isomorphisms. If $P$ is projective both as a left $A$-module
and as a right $B$-module, then the functor $P\otimes_B -$ gives a
Morita equivalence between $B$ and $A$.
\end{prop}
\begin{proof}
We freely use the terminology and results of \S3.5 of \cite{mcr}.
In particular, if $M$ is a left module over a ring $R$, and
$S=\End_{\lMod{R}}(M)^{op}$, then in order to show that $M$ yields
a Morita equivalence between $R$ and $S$, it is enough to prove
that $M$ is a projective $R$-module, and that it is a generator of
the category $\lMod{R}$.

\mbr

Let us write
\[
P^* = \Hom_{\lMod{A}}(P,A) \quad\text{and}\quad
P^\vee=\Hom_{\rMod{B}}(P,B);
\]
both spaces have the natural $(B,A)$-bimodule structures. In view
of the previous paragraph, it suffices to show that $P$ is a
generator of the category $\lMod{A}$, in other words, that the
evaluation map
\begin{equation}\label{e:p-star}
\ev : P\tens_B P^* \rar{} A
\end{equation}
is surjective. The fact that $P$ is a projective right $B$-module
implies that the natural map
\[
\phi : P\tens_B P^\vee \rar{} \End_{\rMod{B}}(P) \xleftarrow{\
\simeq\ } A
\]
given by sending $x\tens f$ to the endomorphism $y\mapsto x\cdot
f(y)$ is surjective. However, it is clear that this map must
factor through the map \eqref{e:p-star}. More precisely, if
$\psi:P^\vee\to P^*$ is defined by $\psi(f)(x)=\phi(x\tens f)$,
then $\phi=\ev\circ(1\circ\psi)$. This shows that \eqref{e:p-star}
is surjective, completing the proof of the proposition.
\end{proof}

\subsection{}\label{ss:uniform} In this subsection we recall,
following \S2.2 of \cite{mcr}, the notion of a uniform module, and
explain how it will be applied in our proof of Theorem
\ref{t:main}. Let $A$ be a ring and $M$ a (left) $A$-module. We
say that a submodule $N\subseteq M$ is {\em essential} if $N\cap
X\neq(0)$ for every nonzero submodule $X\subseteq M$. We say that
$M$ is {\em uniform} if $M\neq(0)$ and every nonzero submodule of
$M$ is essential.
\begin{lem}\label{l:morita-uniform}
Let $B$ be a left Ore domain and $P$ an $(A,B)$-bimodule which
yields a Morita equivalence between $A$ and $B$. Then $P$ is a
uniform $A$-module. In particular, if $P'$ is any torsion-free
nonzero left $A$-module, then every surjective homomorphism
$f:P\to P'$ of $A$-modules must be an isomorphism.
\end{lem}
\begin{proof} Note that the equivalence $\lmod{B}\eqv\lmod{A}$
induced by $P$ maps $B$ to $P$. Also, $B$ is a uniform left
$B$-module by the last remark of \S2.2.5 of loc.~cit. This implies
that $P$ is a uniform left $A$-module, by Lemma 3.5.8(vi) of
loc.~cit. In particular, if $f:P\to P'$ is as in the statement of
the lemma, and $K=\Ker(f)\neq(0)$, then $K$ is an essential
submodule of $P$. If $x\in P$ is any element such that $f(x)\neq
0$, then we see that $ax\in K$ for some $a\in A$, which implies
that $a\cdot f(x)=0$, contradicting the assumption that $P'$ is
torsion-free.
\end{proof}

\subsection{}\label{ss:gorenstein} We now record some observations in noncommutative
ring theory that are interesting in their own right. For the
definition of and basic facts about rigid dualizing complexes, as
well as the notions of Cohen-Macaulay rings and Gorenstein rings
in the noncommutative situation, we refer the reader to \S3 of
\cite{EG} and references therein.
\begin{prop}\label{p:cohen-mac} Let $A$ be a finitely generated connected filtered
algebra such that $\gr A$ is commutative and Gorenstein.
\begin{enumerate}[(1)]
\item If $M$ is a finitely generated left $A$-module equipped with a good
filtration, and $\gr M$ is a \cm $\gr A$-module of degree $k$,
then $M$ is a \cm $A$-module of degree $k$.
\item If, in addition, $A$ is regular and $M$ is a finitely generated \cm
$A$-module of degree $0$, then $M$ is projective.
\end{enumerate}
\end{prop}
\begin{proof}
The argument for (1) is essentially given in \S3 of \cite{EG}.
Namely, since $\gr A$ is commutative and Gorenstein, it is its own
dualizing complex, and then the same is true of $A$. But for every
$j\geq 0$, a result of Bj\"ork (\cite{bjork}, Proposition 3.1)
shows that $\gr\Ext_A^j(M,A)$ is a subquotient of $\Ext_{\gr
A}^j(\gr M, \gr A)$. Thus if $\Ext^j_{\gr A}(\gr M,\gr A)=(0)$ for
all $j\neq k$, then $\Ext_A^j(M,A)=(0)$ for all $j\neq k$, which
proves (1).

\medbreak

For (2), in view of the previous paragraph, it suffices to show
that if $\Ext^j_A(M,A)=(0)$ for all $j\geq 1$, then $M$ is
projective. We use induction on the projective dimension of $M$,
which is finite by assumption. If the projective dimension is $0$,
we are done. Otherwise there is an exact sequence
\begin{equation}\label{e:1}
0 \rar{} N \rar{} P \rar{} M \rar{} 0
\end{equation}
with $P$ finitely generated and projective. The projective
dimension of $N$ is one less than that of $M$, and the long exact
sequence of $\Ext$'s shows that $\Ext_A^j(N,A)=(0)$ for all $j\geq
1$, whence $N$ is projective by induction. In particular, our
assumption on $M$ implies that $\Ext^1_A(M,N)=(0)$, so the
sequence \eqref{e:1} splits, and thus $M$ is projective.
\end{proof}

\subsection{}\label{ss:proof} We return to the situation of
Theorem \ref{t:main}. Given the construction of the minimal
resolution $\widetilde{X}\to X$ discussed in Section
\ref{s:minres}, it is clear that we only need to prove the first
part of the theorem, namely, that the $\bZ$-algebra
$B(\la+\xi,\chi)$ satisfies
\[
\lmod{\cO^\la}\simeq \qgr{B(\la+\xi,\chi)}
\]
for all sufficiently large $\xi$. Using Proposition \ref{p:shifts}
and Theorem \ref{t:morita}, we see that it is enough to check that
$B(\la+\xi,\chi)$ is a Morita $\bZ$-algebra, in the sense of
Definition \ref{d:morita}, for all sufficiently large $\xi$. We
combine the results of Section \ref{s:minres} and
\S\S\ref{ss:o-lambda}--\ref{ss:gorenstein}. The second condition
in the definition of what it means for $B(\la+\xi,\chi)$ to be a
Morita $\bZ$-algebra follows from the first one in view of Theorem
\ref{t:geom}(3) and Lemma \ref{l:morita-uniform}. Thus we need to
check that there exists a large enough $\xi\in\lpp$ for which the
$(\cO^{\la+\xi+i\chi},\cO^{\la+\xi+j\chi})$-bimodules
\[
P_{ij}:=P^{\la+\xi+j\chi}_{(i-j)\chi}
\]
satisfy the assumptions of Proposition \ref{p:morita} whenever
$i-j\gg 0$.

\subsection{}
Let $N\in\bN$ be fixed, so that the statements of the last two
parts of Theorem \ref{t:geom} hold for all $m,n\geq N$. The key
point is that, {\em whereas the bimodules $P_{ij}$ themselves
depend on $\xi$, the associated graded modules do not}. In fact,
Theorem \ref{t:hol2} implies that $\gr P_{ij}\cong S_{i-j}$ as
$R$-bimodules, using the notation of Section \ref{s:minres}. Now
the second statement of Proposition \ref{p:torsfree} and the
standard argument given, e.g., in \S3 of \cite{EG}, in the proof
of Theorem 1.5(iv), implies that for $i-j\geq N$, we have
\[
\cO^{\la+\xi+i\chi} \rar{\simeq}
\End_{\rmod{\cO^{\la+\xi+j\chi}}}(P_{ij})
\]
and
\[
\cO^{\la+\xi+j\chi}\rar{\simeq}
\End_{\lmod{\cO^{\la+\xi+i\chi}}}(P_{ij})^{op}.
\]
Hence we only need to make sure that $P_{ij}$ is projective as a
left and as a right module.

\medbreak

For simplicity, let us write $\cO_i=\cO^{\la+\xi+i\chi}$. Now for
$n=N,N+1,\dotsc,2N-1$, the first statement of Proposition
\ref{p:torsfree} and the argument used in the proof of Proposition
\ref{p:cohen-mac}(1) imply that if $i-j=n$, then the modules
\begin{equation}\label{e:ext1s}
\Ext^\ell_{\lmod{\cO_i}}(P_{ij},\cO_i) \quad\text{and}\quad
\Ext^\ell_{\rmod{\cO_j}}(P_{ij},\cO_j) \qquad (\ell=1,2)
\end{equation}
have finite dimension which is {\em independent of $\xi$}, and
depends only on $i-j$ but not on $i$ or $j$ separately.
Furthermore,
\[
\Ext^\ell_{\lmod{\cO_i}}(P_{ij},\cO_i) =
\Ext^\ell_{\rmod{\cO_j}}(P_{ij},\cO_j) = (0) \qquad \text{for }
\ell\geq 3,
\]
because the global dimensions of $\cO_i$ and $\cO_j$ are at most
$2$ (see \cite{CBH}). In particular, by Proposition \ref{p:findim}
and Remark \ref{r:findimright}, we can choose and fix a large
enough $\xi$ for which the modules \eqref{e:ext1s} are necessarily
zero if $i-j\in\{N,N+1,\dotsc,2N-1\}$. Thus, by Proposition
\ref{p:cohen-mac}(2), we have now shown that, for the $\xi$ that
we've chosen, the bimodules $P_{ij}$ induce Morita equivalences
between the algebras $\cO_j$ and $\cO_i$ whenever
$i-j\in\{N,N+1,\dotsc,2N-1\}$.

\begin{rem}\label{r:xi}
Note that if $\cO^\la$ has no nonzero finite dimensional modules
(by Theorem 0.3 of \cite{CBH}, this happens if and only if
$\la\cdot\al\neq 0$ for all non-Dynkin roots $\al$), then the
modules \eqref{e:ext1s} are automatically zero, and the choice of
$\xi$ is unnecessary. We conjecture that, in fact, the modules
$P_{ij}$ induce Morita equivalences between the algebras $\cO_i$
and $\cO_j$ provided $\la$ is {\em dominant} in the sense that
$\operatorname{Re}(\la\cdot\al)>0$ for every positive Dynkin root $\al$.
\end{rem}

\subsection{} We complete the proof as follows. With the notation
above, it is obvious that any integer $m\geq 2N-1$ can be written
as a sum $m=m_1+\dotsb+m_k$, where each
$m_s\in\{N,N+1,\dotsc,2N-1\}$. Thus, for $i-j=m$, we see from
Theorem \ref{t:geom}(3) that the bimodule $P_{ij}$ is a
homomorphic image of the tensor product
\[
\sP:=P_{j+m_1,j}\otimes P_{j+m_2,j+m_1} \otimes\dotsm\otimes
P_{i,j+m_{k-1}}
\]
over the appropriate algebras $\cO_{j+m_s}$. Since each factor in
the tensor product induces a Morita equivalence by the argument
above, so does the tensor product. Hence, by Lemma
\ref{l:morita-uniform}, the surjection $\sP\twoheadrightarrow
P_{ij}$ must necessarily be an isomorphism; in particular,
$P_{ij}$ also induces a Morita equivalence between $\cO_j$ and
$\cO_i$, and the proof of the theorem is complete.

\end{document}